\documentclass[11pt,epsf]{article}
\usepackage{graphicx}
\usepackage{theorem, amsmath, amssymb}
\newtheorem{Theorem}{Theorem}[section]
\newtheorem{Definition}[Theorem]{Definition}
\newtheorem{Proposition}[Theorem]{Proposition}

\newtheorem{Lemma}[Theorem]{Lemma}

\theoremstyle{remark}
\newtheorem{Remark}[Theorem]{Remark}

  {\begin{equation}\left\{\begin{aligned}}%
  {\end{aligned}\right.\end{equation}\ignorespacesafterend}%

\newenvironment{SDE*}%
  {\begin{equation*}\left\{\begin{aligned}}%
  {\end{aligned}\right.\end{equation*}\ignorespacesafterend}%

\usepackage{color}
\begin{document}
\title{
A note on the  strong formulation of  stochastic control problems\\  with model uncertainty
}

\author{Mihai S\^{\i}rbu \footnote{University of Texas at Austin,
    Department of Mathematics, 1 University Station C1200, Austin, TX,
    78712.  E-mail address: sirbu@math.utexas.edu. The research of
    this author was supported in part by the National Science
    Foundation under Grant    DMS 1211988. Any opinions, findings, and conclusions or recommendations expressed in this material are those of the authors and do not necessarily reflect the views of the National Science Foundation.}}
\maketitle
\begin{abstract} We consider a  Markovian stochastic control problem with  model uncertainty. The controller (intelligent player) observes only the state, and, therefore, uses feed-back (closed-loop) strategies.  The adverse player (nature) who does not have a direct interest in the pay-off, chooses open-loop controls that parametrize Knightian uncertainty. This creates a two-step optimization  problem (like half of a game) over feed-back strategies and open-loop controls. The main result is to show that, under some assumptions, this provides the  same value as the  (half of) the zero-sum symmetric game where the adverse player  also plays feed-back strategies and actively tries to minimize the pay-off. The value function is independent of the filtration accessible to the adverse player.
 Aside from the modeling issue, the present note is a technical companion to \cite{sirbu}.
\end{abstract}
%\vfill


\section{Introduction}
%\subsection{Existing Literature}
We consider a stochastic control  problem with model uncertainty. At first, the  problem looks   identical  to the symmetric zero-sum game in  \cite{sirbu}. However, here, only one player  is a true optimizer (intelligent player) who tries to maximize the pay-off. The other control variable is  chosen by an adverse player (nature) who does not have a vested interest in minimizing the pay-off, and models Knightian uncertainty. We argue that the two apparently identical problems (the symmetric zero-sum game in \cite{sirbu} and the  model uncertainty) should be rigorously defined differently.
%Formally, the robust optimization problem {\bf appears} to be identical to half of the game. When modeling it rigorously, we argue that, this is actually {\bf not} identical to half of a zero-sum game. 

More precisely, we interpret the control problem with model uncertainty as a two-step optimization problem. The controller (intelligent players) observes the state process only, so he/she chooses feed-back (closed-loop) strategies.
The adverse player chooses open-loop controls, and such controls are actually adapted to a possibly larger filtration than the one generated by the Brownian motion. In other words, the adverse player, while not acting strategically against the controller, has access to the Brownian motion and other information and  may choose a parametrization of the model which {\bf just happens} to be  totally adverse to the controller.

A similar model of robust control over feed-back/closed-loop/positional strategies for the controller and open loop-controls for the adverse player has been considered in \cite{krasovskii-subbotin-88} in  deterministic setting. However, our discretization of time for the feed-back strategies is different, and, arguably, better fitted to the present case where the system is stochastic, and allows for strong solutions  of the state system. In addition, our note deals with the (important, in our view) issue of the  information available to the adverse player.  A part of our contribution is to prove that  the value function does not depend on the filtration accessible to the adverse player. This is not obvious a-priori.

There is a vast literature on robust optimization/model uncertainty, and we do not even attempt to scratch the surface in presenting the history of the problem. However, we have not encountered  this very particular way to represent stochastic optimization  problems with  model uncertainty, i.e.   a {\bf strong}  formulation over {\bf elementary feed-back strategies}  for the controller vs. {\bf open-loop controls for the nature}, nor the technical result about the equality of the value functions we obtain.

The message of the present note is two-fold: first, an optimization problem with model uncertainty is {\bf not} the same as a zero-sum game, so it should be modeled differently. 
We propose to use feed-back strategies for the controller and open-loop controls for the adverse player, obtaining a two-step/sup-inf optimization problem over  {\bf strong solutions} of the state system. 
 Second, with this formulation, the value function is, indeed, equal to the (lower) value of the  zero-sum game, where the adverse player is symmetric to the controller and also plays pure feed-back strategies. Beyond the modeling issue, the mathematical statement does not seem obvious, and the proof is based on verification by Stochastic Perron's Method, along the lines of \cite{sirbu}. It is unclear how one could prove directly, using only the probabilistic representation of the value functions, such statement.
\section{Stochastic Control with Model Uncertainty}
\subsection{The Stochastic System}
We consider a stochastic differential system of the form:
\begin{equation}\label{eq:SDE}\left\{
\begin{array}{ll}
 dX_t=b(t,X_t,u_t, v_t)dt+\sigma (t, X_t,u_t, v_t)dW_t,   \\ 
 X_s=x \in \mathbb{R}^d,
 \end{array} \right.
 \end{equation}
 starting at an initial time $0\leq s\leq T$ at some  position $x\in \mathbb{R}^d.$
 Here, the control $u$ chosen by the controller (intelligent player) belongs to some compact metric space $(U,d_U)$ and 
 the parameter $v$ (chosen by the adverse player/nature) belongs to some other  compact metric space $(V,d_V)$ and  represents the model uncertainty. In other words, the Brownian motion $W$ represents the ``known unknowns'', and the process $v$ stands for  the ``unknown unknowns'', a.k.a. ``Knightian uncertainty''. The state $X$ lives in $\mathbb{R}^d$ and  the process $(W_t)_{s\leq t\leq T}$ is a $d'$-dimensional Brownian motion on a \emph{fixed} probability space
 $(\Omega, \mathcal{F}, \mathbb{P})$ with respect to some filtration
 $\mathbb{F}=(\mathcal{F}_t)_{s\leq t\leq T}$. The filtration $\mathbb{F}$ satisfies the usual conditions and is  is usually larger than the  the augmented natural filtration generated by the Brownian motion, by which we mean,
 $\mathcal{F}^W_t=\sigma (W_u, s\leq u\leq t)\vee \mathcal{N}(\mathbb{P}, \mathcal{F})\ \  \textrm{for}\ \ s\leq t\leq T.$
 The space $(\Omega, \mathcal{F}, \mathbb{P})$, the Brownian motion $W$ and the filtration $\mathbb{F}$  may depend on $s$.  To keep the notation simple, we do not emphasize the dependence on $s$, unless needed.
The coefficients $b:[0,T]\times \mathbb{R}^d\times U\times V\rightarrow \mathbb{R}^d$ and 
$\sigma :[0,T]\times \mathbb{R}^d\times U\times V\rightarrow \mathcal{M}^{d\times d'}$ satisfy the 

\noindent {\bf Standing assumption:}
\begin{enumerate}
\item {\bf (C)}  $b, \sigma$  are jointly continuous on $[0,T]\times \mathbb{R}^d\times U\times V$

\item {\bf (L)} $b,\sigma$ satisfy a uniform local Lipschitz condition in $x$, i.e.
\begin{equation*}\label{Lip}
|b(t,x,u,v)-b(t,y,u,v)|+|\sigma (t,x,u,v)-\sigma (t,y,u,v)|\leq L(K) |x-y|\ \ 
\end{equation*}
$\forall\  |x|,|y|\leq K, t\in [0,T],\  u\in U, v\in V$
for some $L(K)<\infty$, and 
\item {\bf (GL)}
$b,\sigma$ satisfy a global linear growth condition in $x$
$$  
|b(t,x,u,v)|+|\sigma(t,x,u,v)|\leq C(1+|x|)$$
$\forall\  |x|,|y|\in \mathbb{R}^d, t\in [0,T],\  u\in U, v\in V$ for some $C<\infty.$
 \end{enumerate}
 Now, given a   bounded and continuous function $g:\mathbb{R}^d\rightarrow \mathbb{R}$, the controller is trying to maximize
 $\mathbb{E}[g(X^{s,x;u,v}_T)].$ Since $v$ is "uncertain", optimizing ``robustly'', means  optimizing the functional $\inf _v\mathbb{E}[g(X^{s,x;u,v}_T)]$, leading to the two-step optimization problem $$\sup _{u }\left (\inf _v\mathbb{E}[g(X^{s,x;u,v}_T)]\right).$$
 It is not yet clear what $u$, $v$ mean in the formulation above, and giving a  precise meaning to this is  one of the goals of the present note.

\subsection{Modeling a Zero-Sum Game} For an identical stochastic system,  imagine that $v$ represents the choice of another  intelligent player and  $g(X^{s,x;u,v}_T)$ is the amount payed by the $v$ player to the $u$ player. For this closely related, but  {\bf different} problem it was argued in \cite{sirbu} that,
{as long as both players only observe the state process, they should both play, symmetrically, as strategies, some feed-back functionals $u,v$ of restricted form.

We denote by $C([s,T])\triangleq C([s,T],\mathbb{R}^d)$ and endow this path space with the natural (and raw) filtration 
 $\mathbb{B}^s=(\mathcal{B}^s_t)_{s\leq t\leq T}$ defined by
 $\mathcal{B}^s_t\triangleq \sigma (y(u),s\leq u\leq t), \ \ s\leq t\leq T.$
 The elements of the path space $C([s,T])$ will be denoted by $y (\cdot)$ or $y$. The stopping times on the space $C([s,T])$ with respect with the filtration $\mathbb{B}^s$, i.e.  mappings $\tau :C([s,T])\rightarrow [s,T]$ satisfying
 $  \{\tau \leq t\}\in \mathcal{B}^s_t \ \forall \ s\leq t\leq T$
are called  stopping  rules, following \cite{ks}. We denote   by $ \mathbb{B}^s$ the class of such stopping  rules starting at $s$.
 \begin{Definition}[Elementary Feed-Back Strategies] \label{def:s} Fix $0\leq s\leq T$.
 An  elementary  strategy $\alpha$ starting at $s$,  for the first intelligent player/controller  is defined by

\begin{itemize}
\item a finite non-decreasing sequence of stopping  rules, i.e. $\tau _k \in \mathbb{B}^s$ for $k=1, \dots, n$ and 
$$s=\tau _0\leq \dots \tau _k\leq \dots \leq \tau _n=T $$
\item for each $k=  1\dots n$, a constant value of the strategy $\xi_k$ in between the times $\tau _{k-1}$ and $\tau _k$, which is decided based only on the knowledge of the past state up to  $\tau _{k-1}$, i.e.
$\xi_k:C([s,T])\rightarrow U$ such that
$\xi _k\in \mathcal{B}^s_{\tau _{k-1}}$.
\end{itemize}
 The strategy is to hold $\xi_k$ in between $(\tau _{k-1}, \tau _{k}]$, i.e.
 $\alpha: (s,T]\times C([s,T])\rightarrow U$ is defined by
 $$\alpha(t, y(\cdot))\triangleq \sum _{k=1}^n \xi _{k}(y(\cdot))1_{\{ \tau _{k-1}(y(\cdot))<t\leq  \tau _{k}(y(\cdot))\}}.$$
  An elementary strategy $\beta$ for the second player is defined in an identical way, but takes values in $V$.
 We denote by $\mathcal{A}(s)$
  and $\mathcal{B}(s)$ 
 the collections of all possible elementary strategies for the  $u$-player and the $v$-player, respectively, given the initial deterministic time $s$.
\end{Definition}
The main result in \cite{sirbu} is the description of the lower and upper values of such a {\bf zero-sum symmetric game over elementary feed-back strategies}. We recall below the result, for convenience:
\begin{Theorem} Under the standing assumption, we have
\begin{enumerate}
\item for each $\alpha\in \mathcal{A}(s)$, $\beta \in \mathcal{B}(s)$, there exists a unique strong solution $(X^{s,x;\alpha ,\beta}_t)_{s\leq t\leq T}$ (such that  $X^{s,x;\alpha,\beta}_t\in \mathcal{F}^W_s$) of the closed-loop state system
\begin{equation}\label{eq:see-closed}\left\{
\begin{array}{ll}
 dX_t=b(t,X_t, \alpha (t, X_{\cdot}), \beta (t,X_{\cdot}))dt+\sigma (t, X_t,\alpha (t, X_{\cdot}), \beta (t,X_{\cdot}))dW_t,   \ s\leq t\leq T\\ 
 X_s=x \in \mathbb{R}^d.
 \end{array} \right.
 \end{equation}
\item the functions 
$$V^-(s,x)\triangleq \sup _{\alpha \in \mathcal{A}(s)}\inf _{\beta \in \mathcal{B}(s)}\mathbb{E}[g(X^{s,x;\alpha,\beta}_T)]\leq 
V^+(s,x)\triangleq \inf _{\beta \in \mathcal{B}(s)}  \sup _{\alpha \in \mathcal{A}(s)}\mathbb{E}[g(X^{s,x;\alpha,\beta}_T)]$$
are the unique bounded continuous viscosity solutions of the  Isaacs equations (for $i=-$ and $i=+$) to the game
\begin{equation}\label{eq:Isaacs}
\left \{
\begin{array}{ll}
-v_t-H ^{i}(t,x,v_x,v_{xx})=0\ \ \textrm{on}\ [0,T)\times \mathbb{R}^d,\\
v(T,\cdot)=g(\cdot),\ \ \textrm{on}\ \mathbb{R}^d.
\end{array}
\right.
\end{equation}
where, \[
H^-(t,x,p,M)\triangleq \sup_{u \in U} \inf _{v\in V} L(t,x,p,M;u,v)\leq 
H^+(t,x,p,M)\triangleq \inf _{v\in V}\sup_{u \in U}L(t,x,p,M;u,v),
\] 
using the notation
$L(t,x,p,M;u,v)\triangleq
b(t,x,u,v)\cdot p+\frac{1}{2}Tr \left (\sigma(t,x,u,v)\sigma(t,x,u,v)^T M\right).$\end{enumerate}
\end{Theorem}
\subsection{ Back to Control with Model Uncertainty}
In our setting, $v$ does not represent an intelligent player: we can think about it as  nature,  which does not have a pay-off to minimize (or  a vested interest from playing against player $u$). The controller (player $u$) does have a pay-off to maximize.  It is still natural to assume that, the controller only observes the state of the system, so he/she uses {\bf the same elementary feedback strategies} $\alpha \in \mathcal{A}(s)$.  On the other hand, the adverse player, the nature, can choose any parameter $v$, and, can actually do so using the whole information available in the filtration $\mathbb{F}$. In other words, we treat as  the possible (uncertain) choices of the model to be all {\bf open-loop} control processes
$v_t$. We define
$$\mathcal{V}(s)\triangleq\{(v_t)_{s\leq t\leq T}|\textrm{predictable with respect to } \mathbb{F} \},$$
and set up the optimization problem under model uncertainty as
$$V(s,x)\triangleq \sup _{\alpha \in \mathcal{A}(s)}\inf _{v\in \mathcal{V}(s)}\mathbb{E}[g(X^{s,x;\alpha ,v}_T)].$$
The above formulation represents the modeling contribution of the present note.
We emphasize  one last time that, in our model,
\begin{itemize}
\item nature uses open-loop controls $v\in \mathcal{V}(s)$, while the controller uses feed-back strategies $\alpha \in \mathcal{A}(s)$,
\item the nature's controls are adapted to the filtration $\mathbb{F}$ which may be strictly larger than the one generated by the Brownian motion.
\end{itemize}
Before even studying the well posed-ness of the state equation over one feed-back strategy $\alpha$ and one open loop control $v$, it is expected (proven rigorously below), that
$V\leq V^-$.
The main result of the present note   is

\begin{Theorem}\label{thm:main} Under the standing assumption, we have
\begin{enumerate}
\item for each $\alpha \in \mathcal{A}(s)$ and $v\in \mathcal{V}(s)$, the state equation has a unique strong solution $(X^{s,x;\alpha,v}_t)_{s\leq t\leq T}$ with  $X^{s,x;\alpha,v}_t\in \mathcal{F}_t\supset \mathcal{F}^W_t,$
\item $V=V^-$ is the unique continuous viscosity solution of the lower Isaacs equation,
\item the value function $V$ satisfies the Dynamic Programming Principle
$$ V(s,x)= \sup_{\alpha \in \mathcal{A}(s)}  \inf _{v \in \mathcal{V}(s)}\mathbb{E}[g(X^{s,x;\alpha, v}_{\rho(X^{s,x;\alpha,v})})]\ \forall \ \rho \in \mathbb{B}^s.$$
\end{enumerate}

\end{Theorem}
It is important, in our view, to obtain strong solutions of the state equation, and this is the main reason to restrict feed-back strategies to the class of elementary strategies.
%In words, the main theorem says that, in our set-up, the nature is the stronger player, despite not having a direct interest in minimizing the pay-off. 
Mathematically, our result states that the use of open-loop controls by the stronger player (here, the nature), even adapted to a much larger filtration than the one generated by the ``known randomness'' $W$, does not change the value function, from the one where the stronger player only observes the state process, as long as the weaker player {\bf only observes the state}.
More precisely, the technical contribution of the note is to show that 
$$ \sup _{\alpha \in \mathcal{A}(s)}\inf _{v\in \mathcal{V}(s)}\mathbb{E}[g(X^{s,x;\alpha,v}_T)]=
 \sup _{\alpha \in \mathcal{A}(s)}\inf _{\beta \in \mathcal{B}(s)}\mathbb{E}[g(X^{s,x;\alpha ,\beta }_T)].$$
 In our understanding, this is not entirely obvious.
% The other part of the contribution is to {\bf modeling the robust control/model uncertainty}  problem as  
%$$\sup _{\alpha\in \mathcal{A}(s)}\left (\inf _{v\in \mathcal{V}(s)}\mathbb{E}[g(X^{s,x;\alpha,v}_T) ]\right ).$$
\begin{Remark}
\begin{enumerate}
\item one possible way to model the robust control  problem is to assume that $\alpha$ is an Elliott-Kalton strategy  (like in \cite{ek} or \cite{fs}) and $v$ is an open loop control. 
%While this would work out mathematically  fine (and can actually be found in the literature)
%\cite{bouchard-nutz}
While such an approach is present in the literature,
 we find it  quite hard to justify the assumption  that the controller can observe the changes in model uncertainty {\bf in real time}, i.e. really observe $v_t$ right at time $t$.  Locally (over an infinitesimal time period), this amounts for the nature to {\bf first} choose the uncertainty parameter $v$, {\bf then, after observing $v$} for the controller to choose $u$.  This contradicts the very  idea of Knightian uncertainty we have in mind.
 % In different setting, such a game is modeled, for example, in \cite{bouchard-nutz}. 
 If one actually  went ahead and modeled our control problem in such a way, than $V$ would be equal to $V^+$, since the Elliott-Kalton player is the stronger player  as described above (see \cite{fs}  for the mathematics, under stronger assumptions on the system). 
 
 \item another way would be to model the ``nature'' as the Elliott-Kalton strategy player $\beta$ an let the controller/intelligent player use open loop controls $u$. This does not seem too appealing either, since nature does not have any pay-off/vested interest. Why would nature be able to {\bf observe the controller's actions and act strategically against} him/her? In addition, if the controller chooses open-loop controls, he/she needs to have the whole information in $\mathbb{F}$ available.  The controller does not  usually observe directly even the noise $W$, leave alone the other possible information in $\mathbb{F}.$
%However,  if  one simply modeled the problem this way,  the controller, would need to chooes open loop controls $u$ which is also hard to justify, since this player usually does not have access to the whole filtration $\mathbb{F}^s$. 
However, with such a model, mathematically, the resulting value function is expected to  be the same, $V=V^-$ (see, again, \cite{fs}, up to technical details).
\end{enumerate}
\end{Remark}

\section{Proofs}
The proposition below contains the proof of the first item in Theorem \ref{thm:main}.

\begin{Proposition}\label{prop:state-eq-simple} Fix $s,x$ and 
 $\alpha \in \mathcal{A}(s)$ and $ v \in\mathcal{V}(s)$. Then, there exists a unique strong (and square integrable) solution
$(X^{s,x;\alpha,v}_t)_{s \leq t\leq T}$,  $X^{s,x;\alpha,v}_t\in \mathcal{F}_t$ of the state equation
\begin{equation}
\left \{
\begin{array}{ll}
dX_t=b(t,X_t,  \alpha(t, X_{\cdot}), v_t)dt+\sigma (t, X_t,\alpha (t, X_{\cdot}), v_t)\, dW_t,\ s \leq t\leq T\\
X_{s}=x\in \mathbb{R}^d.
\end{array}
\right .
\end{equation}
\end{Proposition}
Proof: The proof of the above proposition (both existence and uniqueness) is based  on solving the equation, successively on $[\tau _k \left (X^{s,x;\alpha,v}_{\cdot}\right),  \tau _{k+1}\left (X^{s,x;\alpha,v}_{\cdot}\right)]$ for $k=1, \dots,n$. The details are rather obvious and, even in \cite{sirbu}, the proof of a similar lemma (where both players choose elementary feed-back strategies, unlike here) was only sketched. $\diamond$

 %together with the following very simple  lemma from \cite{sirbu}
%\begin{Lemma}\label{lemma:stopping} Fix $s$ and 
%let $\tau$ be a stopping  rule,
%$\tau :C([s,T])\rightarrow [s,T]$, $\tau \in \mathbb{B}^s$. Let $(X_t)_{s\leq t\leq T}$  be a process with continuous (all, not only almost surely) paths, which is adapted to $\mathbb{F}.$ Then,  the random time
%$\tau _X:\Omega \rightarrow [s,T]$ defined by
%$\tau _X(\omega)\triangleq \tau (X_{\cdot}(\omega))$ is a stopping time w.r.to the filtration $\mathbb{F}.$ In addition$X_{\tau _X}\in \mathcal{F} _{\tau _X}.$\end{Lemma}
Before we proceed, let $\alpha \in \mathcal{A}(s), \beta \in \mathcal{B}(s)$. We can consider
$v_t=\beta (t, X^{s,x;\alpha, \beta}_{\cdot})\in \mathcal{V}(s),$ such that
$$X^{s,x;\alpha, \beta}_{\cdot}=X^{s,x;\alpha, v}_{\cdot}.$$
This means that, for a fixed $\alpha$, there are more open-loop control nature can use, than feed-back strategies an adverse zero-sum player could use.
This shows that
$$V(s,x)= \sup _{\alpha \in \mathcal{A}(s)}\inf _{v\in \mathcal{V}(s)}\mathbb{E}[g(X^{s,x;\alpha,v}_T)]\leq 
 \sup _{\alpha \in \mathcal{A}(s)}\inf _{\beta \in \mathcal{B}(s)}\mathbb{E}[g(X^{s,x;\alpha ,\beta }_T)]=V^-(s,x).$$
The goal is to prove the inequality above is actually a true equality. 
The proof of the main Theorem \ref{thm:main} relies on a similar  adaptation of the  Perron's Method  that  was introduced in \cite{sirbu} for symmetric zero-sum games played over elementary feed-back strategies. As mentioned, the present note  is a technical companion to \cite{sirbu}. The main (but not only) technical difference is that the stochastic sub-solutions of the robust control problem need to be defined  differently, to account for the fact the the adverse player is using open-loop controls.

%Compared to it, there are two  technical differences, that we      point out right away:
%\begin{enumerate}
%\item the stochastic sub-solutions of the robust control problem need to be defined  differently, to account for the fact the the adverse player is using open-loop controls,
%\item once stochastic sub-solutions are defined,  the Perron scheme seems identical to \cite{sirbu}. However,  one has to make sure that aplying It\^o formula still leads to super-martingales  when super-posing test functions over the state process. This difference, however, is less important.
%\end{enumerate}

Following \cite{sirbu}, we first define elementary feed-back strategies starting at  sequel  times to the initial (deterministic) time $s$. The starting time is  a stopping rule.
\begin{Definition}[Elementary Strategies starting later] \label{def:tau}Fix $s$ and let 
 $\tau \in \mathbb{B}^s$ be  a stopping  rule.  An  elementary strategy, denoted by  $\alpha \in \mathcal{A}(s,\tau)$,  for the first player, starting at $\tau$, is defined by

\begin{itemize}
\item (again) a finite non-decreasing sequence of stopping  rules, i.e.
$\tau _k \in \mathbb{B}^s,$ $k=1,\dots n$ 
%C([s,T])\rightarrow [s,T],\ \ \  \{\tau _k\leq t\}\in \mathcal{B}^s_t\triangleq \sigma (x(u),s\leq u\leq t),\ \ \ k=0\dots n$$
for some finite $n$, and with 
$\tau =\tau _0\leq \dots \tau _k\leq \dots \leq \tau _n=T.$
\item for each $k=  1\dots n$, a constant action $\xi_k$ in between the times $\tau _{k-1}$ and $\tau _k$, which is decided based only on the knowledge of the past state up  $\tau _{k-1}$, i.e.
$\xi_k:C([s,T])\rightarrow U$ such that
$\xi _k\in \mathcal{B}^s_{\tau _{k-1}}$.
\end{itemize}
 The strategy is, again, to hold $\xi_k$ in between $(\tau _{k-1}, \tau _{k}]$, i.e.. 
 $$\alpha: \{(t,y)| \tau (y)<t\leq T, y\in  C([s,T])\}\rightarrow U \ \textrm{
 with }\ \ 
 \alpha (t, y(\cdot))\triangleq \sum _{k=1}^n \xi _{k}(y(\cdot))1_{\{ \tau _{k-1}(y(\cdot))<t\leq  \tau _{k}(y(\cdot))\}}.$$ The notation is consistent with $\mathcal{A}(s)=\mathcal{A}(s,s)$.
  \end{Definition}
 We recall, still from \cite{sirbu}, that 
 strategies in $\mathcal{A}(s,\tau)$ cannot be used by themselves for the game starting at $s$, but have to be concatenated with other strategies. %More precisely we have
 \begin{Proposition}[Concatenated elementary feed-back strategies]\label{prop:conc} Fix $s$ and 
 let $\tau \in \mathbb{B}^s$ be  a stopping  rule  and $\tilde{\alpha }\in \mathcal{A}(s,\tau).$ Then, for each $\alpha \in \mathcal{A} (s,s)$, the mapping
 $\alpha  \otimes_{\tau} \tilde{\alpha }:(s,T]\times C([s,T])\rightarrow U $
 defined by 
 $$\big (\alpha  \otimes _{\tau} \tilde{\alpha } \big )(t,y(\cdot))\triangleq \alpha (t, y(\cdot)) 1_{\{s<t\leq \tau (y(\cdot))\}}+\tilde{\alpha }(t, y(\cdot))1_{\{\tau (y(\cdot))<t\leq T\}}$$
 is a simple strategy starting at $s$, i.e. 
 $\alpha  \otimes _{\tau} \tilde{\alpha }\in \mathcal{A}(s,s).$  \end{Proposition}
%  \begin{Proposition}[Concatenated open loop controls]\label{prop:conc} Fix $s$ and 
% let $\tau ' $ be  a stopping  time of $\mathbb{F}^s$ and $v, \tilde{v}\in \mathcal{V}(s)$. Then, the process
 %$v  \otimes_{\tau '} \tilde{v }:(s,T]\times \Omega \rightarrow V $
% defined by 
 %$$\big (v  \otimes _{\tau '} \tilde{v } \big )_t\triangleq v_t 1_{\{s<t\leq \tau '\}}+\tilde{v}_t1_{\{\tau '<t\leq T\}} $$
% is an open-loop control  starting at $s$, i.e. $v  \otimes _{\tau '} \tilde{v } \in \mathcal{V}(s)$. \end{Proposition}
 %The proof is quite easy, once we observe that, both strategies $\alpha

%\subsection{Stochastic Perron's Method for Robust Control}
Compared to \cite{sirbu} the definition below has to be carefully modified.
\begin{Definition}[Stochastic Sub-Solution]\label{sub-solution}
A function $w:[0,T]\times \mathbb{R}^d \rightarrow \mathbb{R}$ is called a stochastic sub-solution if
\begin{enumerate}
\item it is bounded, continuous  and $w(T, \cdot)\leq g(\cdot)$,
\item for each $s$ and  for each stopping  rule $\tau \in \mathbb{B}^s$ 
 there exists   an elementary strategy $\tilde {\alpha }\in \mathcal{A}(s,\tau)$ such that, for any
$\alpha \in \mathcal{A}(s)$,  any $v\in \mathcal{V}(s)$, any $x$ and  each stopping  rule  $\rho \in \mathbb{B}^s$, $\tau\leq \rho\leq T$,
with the simplifying  notation
$X\triangleq X^{s, x,\alpha  \otimes _{\tau} \tilde{\alpha}, v }$ and 
$\tau '\triangleq \tau (X), \rho '\triangleq \rho (X),$
we have 
$$w(\tau ', X_{\tau '})\leq \mathbb{E}[w(\rho', X_{\rho'})|\mathcal{F}_{\tau'}]\ \ \mathbb{P}- a.s.$$
\end{enumerate}\end{Definition}
%It is another easy observation that, a  stochastic sub-solution satisfies
%$w\leq V\leq V^-$. Since we have already characterized $V^-$ as the unique solution of the lower Isaacs equation in \cite{sirbu}, it turns out that, we actually need {\bf only half} of the Perron contraction here. 
%\begin{Definition}[Stochastic Super-Solution]
%A function $w:[0,T]\times \mathbb{R}^d \rightarrow \mathbb{R}$ is called a stochastic super-solution  if
%\begin{enumerate}
%\item it is bounded, continuous  and $w(T, \cdot)\geq g(\cdot)$,
%\item for each $s$ and  for each stopping time  $\tau ' :\Omega \rightarrow [s,T]$,  and each strategy $v\in \mathcal{V}(s,s)$, for each $x$, 
% there exists   an open loop control $\tilde {\bf u}\in \mathbb{V}(s,\tau ')$ (depending on $v$ and $\tau '$ and $x$) such that, for any 
%${\bf v} \in \mathbb{V}(s,s)$  and   as well as each stopping time    $\rho ' $, $\tau\leq \rho\leq T$, with the simplifying  notation
%$X\triangleq X^{s, x,u, {\bf v} \otimes _{\tau} \tilde{{\bf v}} }$ and 
%we have 
%$$w(\tau ', X_{\tau '})\geq \mathbb{E}[w(\rho', X_{\rho'})|\mathcal{F}^s_{\tau'}]\ \ \mathbb{P}- a.s.$$
%\end{enumerate}\end{Definition}
Let $w$  a stochastic sub-solution.
 Fix $s$. There exists $\tilde{\alpha }\in \mathcal{A}(s)$  such that,  for each $x$, each $\rho \in \mathbb{B}^s$  and each $v\in \mathcal{V}(s)$ we have 
\begin{equation}\label{half-dpp-w}w(s,x)\leq \mathbb{E}\left [w(\rho (X^{s,x,\tilde{\alpha }, v}_{\cdot}), X^{s,x, \tilde{\alpha },  v}_{\rho (X^{s,x,\tilde{\alpha }, v}_{\cdot})})|\mathcal{F}_s   \right], \ \mathbb{P}-a.s.
\end{equation}
Taking the expectation it is obvious that, if $w$ is a stochastic super-solution, then  we have the half DPP/sub-optimality principle 
\begin{equation}\label{dpp-upp-low}
w(s,x)\leq  \sup _{\alpha \in \mathcal{A}(s,s)} \inf _{v\in \mathcal{V}(s)}\mathbb{E}\left [w(\rho(X^{s,x,\alpha,v}_{\cdot}), X^{s,x,\alpha ,v}_{ \rho (X^{s,x,\alpha,v}_{\cdot}   )})\right ], \ \ \ \forall \rho \in \mathbb{B}^s.
\end{equation} Since $w(T, \cdot)\leq g(\cdot)$, we obtain
$w(s,x)\leq V(s,x)\leq V^-(s,x).$

We have already characterized $V^-$ as the unique solution of the lower Isaacs equation in \cite{sirbu}. Therefore, we actually need {\bf only half} of the Perron construction here. 
We denote by $\mathcal{L}$ the set of stochastic sub-solutions in Definition \ref{sub-solution} (non-empty from the boundedness assumptions). Define
%$$v^- \triangleq \sup _{w\in \mathcal{U}^-} w\leq  V^+\leq \inf _{w\in \mathcal{U}^+}w\triangleq v^+,$$
$$w^- \triangleq \sup _{w\in \mathcal{L}} w\leq  V\leq V^- .$$

\begin{Proposition}[Stochastic Perron for Robust Control]\label{prop:main} Under the standing assumptions, $w^-$ is a LSC viscosity super-solution of the lower Isaacs equation, up to $t=0$.
\end{Proposition}
The following lemmas are very similar  to their counterparts in \cite{sirbu}.
\begin{Lemma}
If $w_1, w_2\in \mathcal{L}$ then $w_1\vee w_2\in \mathcal{L}$.
\end{Lemma}
Fix $\tau \in \mathbb{B}^s$  a stopping rule. Let $\tilde{\alpha}_1, \tilde{\alpha}_2 \in \mathcal{A}(s, \tau)$  be the two  feed-back strategies of the controller, starting at $\tau$ corresponding the the sub-solutions $w_1$ and $w_2$ for  the Definition \ref{sub-solution}. The new strategy starting at $\tau$ defined by 
$$
\tilde{\alpha}(t, y(\cdot))=\tilde{\alpha }_1 (t, y (\cdot)) \, 1_{\{w_1(\tau (y), y (\tau (y)))\geq w_2(\tau (y), y (\tau (y)))\}}+\tilde{\alpha}_2  (t, y(\cdot))\, 1_{\{w_1(\tau (y), y (\tau (y)))< w_2(\tau (y), y (\tau (y)))\}}$$
does the job for the definition of $w\triangleq w_1\vee w_2$ as a stochastic sub-solution .  $\diamond$
\begin{Lemma}There exists a non-decreasing sequence $\mathcal{L}\ni w_n\nearrow w^-$.\end{Lemma}
Proof: according to Proposition 4.1 in \cite{bs-1}, there exist $\tilde{w}_n\in \mathcal{L}$ such that
$w^-=\sup _n \tilde{w}_n.$ Now, we can just define
$w_n =\tilde{w}_1\vee\dots \vee \tilde{w}_n\in \mathcal{L}\nearrow w^-.$  $\diamond$

\noindent {\bf Proof of Proposition \ref{prop:main}} The proof is  similar to \cite{sirbu}. Since It\^o formula applies the same, {\bf regardless of filtration}, it produces sub-martingales in a similar way, even though $v$ is an open-loop control, and the filtration may be larger than the one generated by $W$. This is the key point that allows us to obtain the result.
We only sketch some key points of   the proof, in order to avoid repeating all the similar arguments  in \cite{sirbu}.

 The interior super-solution property for $w^-$: Let $(t_0,x_0)$ in the parabolic  interior $[0,T)\times \mathbb{R}^d$ such that a smooth function $\varphi$ strictly touches $v^+$ from below at $(t_0,x_0)$. Assume, by contradiction, that
$\varphi _t +H^-(t,x,\varphi _x, \varphi _{xx})>0\ \  \textrm{ at}\ \  (t_0, x_0).$
In particular, there exists $\hat{u}\in U$ and $\varepsilon >0$ such that
$$\varphi _t (t_0,x_0)+ \inf_{v\in V}\left[b(t_0, x_0,\hat{u}, v)\cdot \varphi _x(t_0,x_0)+\frac{1}{2}Tr(\sigma(t_0,x_0,\hat{u},v)\sigma(t,x,\hat{u},v)^T \varphi _{xx}(t_0,x_0)) \right] > \varepsilon.$$
To simplify notation, all small balls here are actually included in (i.e. intersected with)   the parabolic interior.
 Since $b, \sigma$  are continuous, and $V$ is compact,  the uniform continuity of the above expression in $(t,x,v)$ for $(t,x)$ around $(t_0, x_0)$ implies that  there exists a smaller $\varepsilon >0$ such that
$$\varphi _t(t,x)+ \inf_{v\in V}\left[b(t, x,\hat{u}, v)\cdot \varphi _x(t,x)+\frac{1}{2}Tr(\sigma(t,x,\hat{u},v)\sigma(t,x,\hat{u},v)^T \varphi _{xx}(t,x)) \right] >\varepsilon,\ \ {\textrm on}\ B(t_0, x_0, \varepsilon).$$
Now, on the compact (rectangular) torus $\mathbb{T}= \overline{B(t_0, x_0, \varepsilon)}- B(t_0, x_0, \varepsilon/2)$ we have that $\varphi <w^-$ and the max of $\varphi -w^-$ is attained, therefore it is strictly negative. In other words
$\varphi <w^- -\eta$ on $\mathbb{T}$ for some $\eta >0$. Since $w_n\nearrow w^-$, a Dini type argument similar to \cite{bs-2} and \cite{bs-3} shows that, for $n$ large enough we have $\varphi <w_n-\eta/2$. For simplicity, fix such an $n$ and call $v=w_n$. Now, define, for small $\delta <<\eta /2$
$$v^{\delta}\triangleq
\left \{ 
\begin{array}{ll}
(\varphi +\delta)\vee v\ \ {\textrm on}\ \  B(t_0, x_0, \varepsilon),\\
v \ \ \textrm{outside}\ \  B(t_0, x_0, \varepsilon).
\end{array}
\right.
$$
Since $v^{\delta}(t_0,x_0)>w^-(t_0,x_0)$, we have a contradiction if $v^{\delta}\in \mathcal{L}$. 
Fix $s$ and let $\tau\in \mathbb{B}^s $ be  a stopping rule for the  initial time $s$. We need to construct an elementary  strategy $\tilde{\alpha}\in \mathcal{A}(s,\tau)$ in the Definition \ref{sub-solution} of stochastic sub-solution for $w^{\delta}$. We do that  as follows:
since $v$ is a stochastic sub-solution, there exists  an elementary strategy $\tilde{\alpha}_1$ for $v$ starting at $\tau$ that does the job in Definition \ref{sub-solution}. Next, 
\begin{enumerate}
\item if $(\varphi +\delta) >v$ at $\tau$, follow the constant action $\hat{u}.$
\item if $(\varphi +\delta) \leq v$ at $(\tau, X _{\tau})$ follow the strategy $\tilde{\alpha}_1$
\item follow the strategy  defined in 1-2 until the first time $\tau _1$ when $(t,X_t)\in \partial B(t_0, x_0, \varepsilon/2)$. On this boundary,  we know that $v^{\delta}=v$.
\item after this, follow the strategy $\tilde{\alpha}_3\in \mathcal{A}(s, \tau _1)$ corresponding to the stochastic sub-solution $v$ with starting stopping rule $\tau _1$
\end{enumerate}
We follow the same arguments as in \cite{sirbu} to make the above ideas rigorous.  We do obtain 
 $v^{\delta}\in \mathcal{L}$, so we  reached a contradiction. 
The terminal condition property for $w^-$ is proved very similarly. $\diamond$

\noindent {\bf Proof of Theorem \ref{thm:main}:} Recall that the first part was proved by Proposition \ref{prop:state-eq-simple}.

Next, the proof of the second item  is finished, once we use the comparison result from Lemma 4.1 in  \cite{sirbu}.  More precisely, we know that 
$w^-\leq V\leq V^-$ and $w^-$ is a viscosity super-solution and $V^-$ is a viscosity solution of the lower Isaacs equation (from  \cite{sirbu}[Theorem 4.1]). Therefore, according to  \cite{sirbu}[Lemma 4.1], we also have $V^-\leq w^-$, so  $V=V^-$ is the unique viscosity solution. 

Finally, the DPP in Item 3 of Theorem \ref{thm:main} is actually an easy observation based on the fact that the value function $V^-$ satisfies a similar (but not identical DPP), that $V=V^-$ and the half DPP \eqref{half-dpp-w}.$\diamond$
\section{Additional Modeling Comments}
In our (strong) model of robust control, the value function of the intelligent player turns out to be $V=V^-.$ Obviously, one can ask the question: should  this player try to randomize feed-back strategies somehow, to get the potentially better value $V^{mix}$ of the value over mixed strategies (for both players) obtained in \cite{sirbu-2} (but in a  martingale symmetric formulation)?

Modeling mixed feed-back strategies for the controller, and open loop-strategies controls for the adverse player is a highly non-trivial issue, and  not obviously possible in strong formulation (see \cite{sirbu-2} for some comments along these lines, for the case of a zero-sum symmetric game).  In our formulation of optimization with model uncertainty, the maximizing player has to settle with the value $V=V^-$.
However, the controller couldn't do better anyway in one of the  two situations:
\begin{enumerate}
\item when the Isaacs condition over pure strategies is satisfied, i.e.
$$\sup_{u \in U} \inf _{v\in V} L(t,x,p,M; u, v)= \inf _{v\in V} \sup_{u \in U} L(t,x,p,M;u,v)$$ so $V^-=V^{mix}=V^+$
\item in any additional situation when $V^-=V^{mix}<V^+$, i.e. all situations in which (even at the formal level) potential randomization for the $u$ player does not change the Hamiltonian. More precisely, if 
$$\sup_{u \in U} \inf _{\nu \in \mathcal{P}(V)} \int L(t,x,p,M, u, v)\nu (dv)
=\inf _{\nu \in \mathcal{P}(V)} \sup_{u \in U}\int L(t,x,p,M, u, v)\nu (dv),$$ 
since 
$$\sup_{u \in U} \inf _{v\in V} L(t,x,p,M; u, v)=\sup_{u \in U} \inf _{\nu \in \mathcal{P}(V)} \int L(t,x,p,M, u, v)\nu (dv)$$
and 
$$\inf _{\nu \in \mathcal{P}(V)} \sup_{u \in U}\int L(t,x,p,M, u, v)\nu (dv)= \inf _{\nu \in \mathcal{P}(V)} \sup_{\mu \in \mathcal{P}(U)}\int L(t,x,p,M, u, v)\mu (dv)\nu (dv)$$
we have 
$$H^-=H^{mix}\leq H^+$$
although the Isaacs condition over pure strategies may not be satisfied ($H^-<H^+$).  In such a situation, the robust controller cannot expect to get  a better value then $V=V^-$.  A sufficient condition for this is for  the map 
$$u\rightarrow L(t,x,p,M;u,v)$$
to be concave. Up to different modeling of strategies, this is exactly the case in the interesting recent contribution \cite{ttu-fs}.

\end{enumerate}
%\pagebreaka 
 %\thispagestyle{empty}
\bibliographystyle{amsalpha} %The style you want to use for references. 
%\bibliography{references} %The files containing all the articles and books you ever referenced. 
\providecommand{\bysame}{\leavevmode\hbox to3em{\hrulefill}\thinspace}
\providecommand{\MR}{\relax\ifhmode\unskip\space\fi MR }
% \MRhref is called by the amsart/book/proc definition of \MR.
\providecommand{\MRhref}[2]{%
  \href{http://www.ams.org/mathscinet-getitem?mr=#1}{#2}
}
\providecommand{\href}[2]{#2}

\end{document}